\documentclass[final]{siamltex}
\usepackage{amssymb,amsfonts,amsmath,epsfig,amsbsy,pifont}

\newcommand{\LRARR}[4]{{\mbox{ \raise 0.4 mm \hbox{$#1$}}} \;
  \mathop{\stackrel{\displaystyle\longrightarrow}\longleftarrow}^{#3}_{#4}
  \; {\mbox{\raise 0.4 mm\hbox{$#2$}}}}
\newcommand{\tick}{\ding{52}}
\newcommand{\cross}{\ding{54}}
\newcommand{\pFPE}{P_{\rm{FPE}}}
\newcommand{\bS}{{\mathbf S}}
\newcommand{\bV}{{\mathbf V}}
\newcommand{\bW}{{\mathbf W}}
\newcommand{\qlam}{\lambda_{\rm \scriptscriptstyle QSSA}}

\title{Error Analysis of Diffusion Approximation 
Methods for Multiscale Systems in Reaction 
Kinetics\thanks{Submitted to the journal's Computational Methods 
in Science and Engineering section \today.
The research leading to these results has
received funding from the European Research Council under the {\it
European Community}'s Seventh Framework Programme ({\it
FP7/2007-2013})/ ERC {\it grant agreement} No. 239870.}}
  
\author{Simon L. Cotter\thanks{School of Mathematics,
University of Manchester, Oxford Road,
Manchester, M13 9PL, United Kingdom;
e-mail: simon.cotter@manchester.ac.uk.
Simon Cotter was also funded by First Grant Award EP/L023989/1 from EPSRC. }
\and
Radek Erban\thanks{Mathematical Institute, University of Oxford,
Andrew Wiles Building, Radcliffe Observatory Quarter, Woodstock 
Road, Oxford OX2 6GG, United Kingdom; 
e-mail: erban@maths.ox.ac.uk}. Radek Erban's work was also 
supported by a Royal Society University Research Fellowship; 
by a Nicholas Kurti Junior Fellowship from Brasenose College, 
University of Oxford; and by the Philip Leverhulme Prize 
from the Leverhulme Trust.
}

\begin{document}

\maketitle

\begin{abstract} 
Several different methods exist for efficient approximation of paths
in multiscale stochastic chemical systems. Another approach is to use bursts of
stochastic simulation to estimate the parameters of a stochastic
differential equation approximation of the paths. In this paper,
multiscale methods for
approximating paths are used to formulate different strategies for
estimating the dynamics by diffusion processes. We then analyse
how efficient and accurate these methods are in a range of
different scenarios, and compare their respective advantages and
disadvantages to other methods proposed to analyse
multiscale chemical networks.
\end{abstract}

\section{Introduction}\label{sec:intro}
A well-mixed chemically reacting system in a container of volume $V$ 
is described, at time $t$, by its $N$-dimensional state vector
\begin{equation}
{\mathbf X}(t) \equiv [X_1(t),X_2(t),...,X_N(t)],
\label{statevector}
\end{equation}
where $N$ is the number of chemical species in the system and 
$X_i(t) \in {\mathbb N}_0$, $i=1,2,\dots,N$, 
is the number of molecules of the $i$-th chemical species at time 
$t$. Assuming that the chemical system is subject to
$M$ chemical reactions
\begin{equation}
\sum_{i=1}^N \nu^-_{j,i}X_i \,
\overset{k_j}{\longrightarrow}
\, 
\sum_{i=1}^N \nu^+_{j,i}X_i,
\quad j=1,\dots,M,
\label{chemicalsystem}
\end{equation}
the time evolution of the state vector ${\mathbf X}(t)$ can be simulated
by the Gillespie SSA~\cite{Gillespie:1977:ESS} which is described
in Table \ref{SSAtable}. Here, $\nu^+_{j,i}$ and $\nu^-_{j,i}$ are 
stoichiometric coefficients and $\nu_{j,i} = \nu^+_{j,i} - \nu^-_{j,i}.$
Step [{\bf 1}] of the algorithm in Table \ref{SSAtable} requires to
specify propensity functions which are, for mass-action
reaction kinetics, given by
\begin{equation}   
\alpha_j({\mathbf x}) 
= 
k_j 
\exp \left[ \left(
1 - \sum_{i=1}^N \nu_{j,i}^-
\right)
\log V \right] 
\prod_{i=1}^N 
(\nu_{j,i}^-) !
\binom{x_i}{\nu_{j,i}^-},
\quad j=1,2,\dots,M.
\label{propensityfunctions}
\end{equation}

\begin{table}
\framebox{%
\hsize=0.97\hsize
\vbox{
\leftskip 6mm
\parindent -6mm
{\bf [1]} Calculate propensity functions $\alpha_k({\mathbf X}(t))$,
$k=1,2,\dots,M$.

\smallskip

{\bf [2]} Waiting time $\tau$ till next reaction is given by \eqref{eq:tau}.

\smallskip

{\bf [3]} Choose one $j\in\{1,2,\ldots,M\}$, with probability 
$\alpha_j({\mathbf X}(t))/\alpha_0({\mathbf X}(t))$,
and perform reaction $R_j$, by adding $\nu_{j,i}$ to each
$X_i(t)$ for all $i = 1,2,\dots, N$.

\smallskip

{\bf [4]} Continue with step [{\bf 1}] with time $t=t+\tau$.

\par \vskip 0.8mm}
}\vskip 1mm
\caption{{\it The pseudo code for the Gillespie SSA.}
\label{SSAtable}}
\end{table}

\noindent
Given the values of the propensity functions,
the waiting time to the next reaction is given by:
\begin{equation}\label{eq:tau}
  \tau = -\frac{\log\left(u\right )}{\alpha_0({\mathbf X}(t))},
  \quad \mbox{where} \quad \alpha_0({\mathbf X}(t)) 
  = \sum_{k=1}^M \alpha_k({\mathbf X}(t))
\end{equation}
and $u\sim U([0,1])$. The Gillespie SSA in Table \ref{SSAtable}
is an exact method, in that the trajectories simulated using 
this algorithm evolve exactly as described by the corresponding
chemical master equation\footnote{The CME is a high/infinite 
dimensional set of ordinary differential equations which 
describe the time evolution of the probability of being in 
a particular state of the system.} (CME). Equivalent and more
efficient formulations of the Gillespie SSA have been
developed in the literature~\cite{Cao:2004:EFS,Gibson:2000:EES}.
However, in certain circumstances they can still be very 
computationally intensive. For instance, if the system that is being 
simulated has some reactions which are likely to occur many times on 
a timescale for which others are unlikely to happen at all, then 
a lot of computational effort is spent on simulating the fast reactions, 
when a modeller may well be more interested in the results of the 
slow reactions~\cite{Haseltine:2002:ASC}.  
In this paper, we will focus on approximate
algorithms for such fast-slow systems.

We refer to reactions which have high average propensities, and whose 
reactions may occur many times on a time scale for which others are 
unlikely to happen at all, as fast reactions. Slow reactions are those 
reactions which are not fast reactions. In reality, there may be several 
different timescales present in the reactions of a particular system, 
but for simplicity we assume there is a simple dichotomy~\cite{Liu:2012:HMS}. We may be 
interested in analysing the dynamics of the ``slow variable(s)'', which 
are chemical species (or linear combinations of the species) which are 
invariant to the fast reactions, and therefore are changing on a slower 
timescale~\cite{Cotter:2011:CAM}.

Efforts have been made to accelerate the Gillespie SSA for multiscale
systems. The Nested Stochastic Simulation Algorithm (NSSA) is such a
method~\cite{E:2005:NSS}. The reactions are split into ``fast'' and
``slow'' reactions. The idea of the NSSA is to approximate the
effective propensities of the slow reactions in order to compute
trajectories only on the slow state space. This is done by using short
bursts of stochastic simulation of the fast reaction subsystem. The
Slow-Scale Stochastic Simulation Algorithm (SSSSA)~\cite{Cao:2005:SSS}
comes from a similar philosophy. Instead of using stochastic
simulations to estimate the effective propensities of the slow
reactions, they are instead found by solving the CME for
the fast reactions (whilst ignoring the slow reactions). This has the
advantage that it does not require any Monte Carlo integration,
however it is limited to those systems for which the CME
can be solved or well approximated for the fast subsystem, which may
not be applicable to some complex biologically relevant systems.

Both of these methods use a quasi-steady state approximation (QSSA) in
order to speed up the simulation of a trajectory on the slow state
space~\cite{Rao:2003:SCK}. Another approach is to approximate the 
dynamics of the slow variable by a stochastic differential equation (SDE).
One can either use short bursts of the Gillespie SSA on a
range of points on the slow state space to approximate the effective
drift and diffusion of the slow variable~\cite{Erban:2006:EFC}
or the Constrained Multiscale Algorithm (CMA)~\cite{Cotter:2011:CAM} 
which utilises a modified SSA that constrains the trajectories 
it computes to a particular point on the slow state space. 
These algorithms can be further extended to automatic detection
of slow variables~\cite{Erban:2007:VFE,Singer:2009:PRD,Cucuringu:2014:DSV},
but, in this paper, we assume that the division of state space into 
slow and fast variables is a priori known and fixed during the 
whole simulation.

The advantage of the SDE approximation methods~\cite{Cotter:2011:CAM,Erban:2006:EFC}, is that the
estimation of the drift and diffusion terms can be easily
parallelised, giving each process a subset of the grid points on the
slow state space. This means that if a user has access to high
performance computing facilities, then the analysis of a given
system can be computed relatively quickly. This is not the case for
trajectory-based methods. One could run many trajectories in
parallel~\cite{Klingbeil:2011:SPS}, however if the aim is to analyse
slow behaviours such as rare switches between stable regions, each
trajectory will still have to be simulated for a long time before such
a switch is possible, regardless of the number of trajectories being
computed simultaneously.

In this paper, we take the approach that we would like to approximate
the dynamics of the slow variable ${\mathbf S}$ by a continuous SDE, 
in the same vein as other previous
works~\cite{Cotter:2011:CAM,Erban:2006:GRN}. We wish to estimate the
effective drift $\bV$ and diffusion matrix $D$ of the slow variable,
resulting in approximate dynamics given by:
\begin{equation}
\mbox{d}\bS = \bV(\bS) \, \mbox{d}t + \sqrt{2D(\bS)} \, \mbox{d}\bW.
\label{eq:SDE}
\end{equation}
Here $\mbox{d}\bW$ denotes standard canonical Brownian motion in $N$
dimensions. In this paper we will focus on examples where the slow
variable is one dimensional, although these results can be extended
to higher dimensions~\cite{Cotter:2011:CAM}. The one-dimensional
Fokker-Planck equation (FPE) corresponding to SDE (\ref{eq:SDE})
is given by:
\begin{equation}
\frac{\partial p}{\partial t} = \frac{\partial}{\partial s} \left
  (\frac{\partial}{\partial s} \big [D(s)p(s,t) \big ] - V(s)p(s,t) \right ).\label{eq:FPE}
\end{equation}

In Section \ref{sec:exist} we introduce five methods for simulation
of multiscale stochastic chemical systems, including two novel
approaches: the Nested Multiscale Algorithm (NMA) in Section
\ref{sec:NMA} and the Quasi-Steady State Multiscale Algorithm (QSSMA)
in Section \ref{sec:QSSMA}. In Section \ref{sec:LIN} we compare
the efficiency and accuracy of the CMA, NMA and QSSMA for a simple
linear system, for which we have an analytical solution for the
CME. Then in Section \ref{bistsection} we apply  CMA, NMA and QSSMA
to a bimodal example. Finally, we discuss the relative accuracy and 
efficiency of these methods against others proposed in the literature
and summarise our conclusions in Section \ref{sec:conc}.

\section{Multiscale Algorithms}\label{sec:exist}
We review three algorithms (CMA, NSSA and SSSSA) previously studied in
the literature in Sections \ref{sec:CMA}-\ref{sec:SSSSA}. In Section
\ref{sec:NMA} we introduce the NMA, which combines ideas from the CMA
and the NSSA. The QSSMA is then introduced in Section \ref{sec:QSSMA}.
\subsection{The Constrained Multiscale Algorithm}\label{sec:CMA}
The CMA is a numerical method for the
approximation of the slow dynamics of a well-mixed chemical system by
a SDE, of the form
\eqref{eq:SDE} which is for a one-dimensional slow variable $S$
written as follows:
\begin{equation}
\mbox{d}S = V(S) \, \mbox{d}t + \sqrt{2D(S)} \, \mbox{d}W.
\label{eq:SDE1d}
\end{equation}

 The effective drift $V$ and effective diffusion $D$ at
a given point $S=s$ on the slow manifold are estimated using a short
stochastic simulation. This simulation (called the Constrained SSA, CSSA) 
is similar to that seen in the Gillespie SSA for the full system,
although it is constrained to a particular value of the slow
variable $S$. The CSSA is given in Table \ref{CSSA table}
where ${\bf F}$ is the vector of fast variables and $S$ is
the slow variable.
\begin{table}
\framebox{%
\hsize=0.97\hsize
\vbox{
\leftskip 6mm
\parindent -6mm
{\bf [1]} Calculate propensity functions $\alpha_i(t)$,
$i=1,2,\dots,M$.

\smallskip

{\bf [2]} Next reaction time is given by \eqref{eq:tau}.

\smallskip

{\bf [3]} Choose one $j\in\{1,\ldots,M\}$, with probability 
$\alpha_j/\alpha_0$, and perform reaction $R_j$.

\smallskip

{\bf [4]} If $S \ne s$ due to reaction $j$ occurring, then
reset $S=s$ while not changing the value of ${\bf F}$.

\smallskip

{\bf [5]} If $X_i<0$ for any $i$,
 then  revert to the state of the system before the reaction $j$ occurred.

\smallskip

{\bf [6]} Continue with step [{\bf 1}] with time $t=t+\tau$.

\par \vskip 0.8mm}
}
\caption{{\it The Conditional Stochastic Simulation Algorithm
    (CSSA). Simulation starts with $S=s$ where $s$ is a given value of
  the slow variable.}\label{CSSA table}}
\end{table}
To estimate the effective drift and diffusion, statistics are
collected about the type and frequency of the changes $dS$ of the slow
variable which is reset in step {\bf[4]} of the CSSA. 
For a simulation of length $T(s)$, the estimations are
given by
\begin{eqnarray}
  V(s) &=& \frac{1}{T(s)} \sum_{m=1}^{Q(T(s)))}dS_m,\label{eq:V}\\
  D(s) &=& \frac{1}{2T(s)}\sum_{m=1}^{Q(T(s))} \left( dS_m \right)^2,
  \label{eq:D}
\end{eqnarray}
where $dS_m$ is the change in $S$ due to the $m$-th
iteration of the CSSA before the reset is made in step [{\bf 4}],
$T(s)$ is the chosen length of CSSA simulation, and $Q(T(s))$ 
is the number of iterations of the CSSA that are made up to
time $T(s)$.

By computing these quantities over a range of values of the slow
variable, approximations can then be found, using standard methods,
to the solution of the steady-state Fokker Planck
equation \eqref{eq:FPE} for the SDE with drift $V$ and diffusion $D$.

\subsection{The Nested Stochastic Simulation
  Algorithm}\label{sec:NSSA}
The NSSA~\cite{E:2005:NSS} is a method for the
reduction of computational complexity of simulating the slow dynamics
of a multiscale chemical system. The reactions in the system are
partitioned into fast reactions $\{R_{f_1}, R_{f_2}, \ldots R_{f_n}\}$
and slow reactions $\{R_{s_1}, R_{s_2}, \ldots R_{s_m}\}$, where $M =
n+m$ is the total number of different reactions. We are interested in
the occurrences of the slow reactions, but the computational effort is
dominated by the simulation of the fast reactions in the standard
SSA given in Table \ref{SSAtable}. However, some/all of the slow reactions are dependent on the
value of the fast variables. In the NSSA, the effective propensities
of the slow reactions are estimated by using short simulations of only
the fast reactions. Using these effective propensities, the Gillespie
algorithm can be applied to just the slow reactions. In systems for
which the QSSA is reasonable, this
algorithm can simulate trajectories of the slow variables at a
fraction of the computational cost of the full Gillespie SSA. It
effectively reduces the system to a lower dimensional chemical system
where all of the reactions are ``slow'', with reaction rates estimated
(where required) using relatively short 
bursts of stochastic simulation of the fast
reactions from the full system.

\subsection{The Slow-Scale Stochastic Simulation Algorithm}\label{sec:SSSSA}
The SSSSA~\cite{Cao:2005:SSS} similarly aims to
reduce the full system to a system which contains only the slow
reactions. In this algorithm, the effective propensities are
calculated not by stochastic simulation, as in the NSSA, but through
application of the QSSA. For some classes
of fast sub-systems, the effective propensity can be explicitly
calculated. For others, the value can be approximated using formulae
given in~\cite{Cao:2005:SSS}. Since there is no requirement to simulate
the fast sub-system as in the NSSA, the speed-up in simulation of
trajectories as compared with the Gillespie algorithm is very
large. In some cases non-linear equations may need to be solved to
find the first or second moments of the value of the fast quantities,
using Newton's method, but even in the worst case scenario the overall
computational cost is far less than Monte Carlo sampling.

\subsection{The Nested Multiscale Algorithm}\label{sec:NMA}
The NMA is a new method, which allows for
efficient approximation of multiscale systems by a SDE. 
As in the CMA, our aim is to approximate
values for the effective drift and diffusion of the slow variables
within the system on a set of grid points. At each grid point, we
simulate the fast sub-system, which allows us to approximate the
effective propensities for the slow reaction. The drift and diffusion
terms are then given by the chemical FPE~\cite{Cotter:2013:AFE,Erban:2009:ASC,Liao:2014:TMPE} for the
system with only the slow reactions present, with the values for the
effective propensities substituted in.

For example, say we have $n$ slow reactions, with effective
propensities $\{\bar{\alpha}_i\}_{i=1}^n$, and with stoichiometric
vectors $\nu_{i,S}$. Here $\nu_{i,S}$ is the change in the slow
variable due to the reaction $R_i$. In the case that the slow variable
is the $j$th chemical species, then $\nu_{i,j}$, but it may be more
complex if the slow variable depends on several chemical species. Then, for a 1-dimensional slow variable, the drift $V$
and diffusion $D$ are given by:
\begin{eqnarray}
V(s) &=& \sum_{i=1}^n \bar{\alpha}_i(s) \nu_{i,S},\label{effDrift}\\
D(s) &=& \frac{1}{2} \sum_{i=1}^n \bar{\alpha}_i(s) \nu_{i,S}^2. 
\label{effDiff}
\end{eqnarray}

The NMA has the advantage over the CMA that it converges on the
timescale of the fast variables, whereas the CMA converges on the
timescale of the slow variables.

\subsection{The Quasi-Steady State Multiscale Algorithm}\label{sec:QSSMA}
The QSSMA follows a similar theme as the NMA, except this time 
we use the methodology of~\cite{Cao:2005:SSS} to approximate the effective propensity
functions. A QSSA is used to derive the value
of these functions, as in the slow-scale SSA (as detailed in Section
\ref{sec:SSSSA}). Once the effective propensities have been calculated,
the formulae (\ref{effDrift})--(\ref{effDiff}) for the drift $V$ and
diffusion $D$ can once again be used to approximate the dynamics of
the slow variable by a SDE of the form (\ref{eq:SDE1d}).
The QSSMA does not require any stochastic simulation in order to
estimate the effective drift and diffusion functions, and thus we see
remarkable speed-ups when compared with the CMA. However, as with the
NMA, other errors come into play that are not present in the
approximations arising from the CMA.

\section{Efficiency and Accuracy of the Schemes}\label{sec:LIN}
In this section, we aim to test the efficiency and accuracy of the
three schemes (the CMA, and the newly proposed NMA and SSMA). To test
the algorithms, we choose a simple multiscale system of two chemical
species $X_1$ and $X_2$ in a reactor of volume $V$ undergoing the following four reactions: 
\begin{eqnarray}
R_1 &:& \qquad\qquad \emptyset \, \overset{k_1}{\longrightarrow} \, X_1 
\nonumber \\
R_2 &:& \qquad\qquad X_2 \, \overset{k_2}{\longrightarrow} \, \emptyset 
\label{eq:lin} \\
R_3 &:& \qquad\qquad X_1 \, \overset{K}{\longrightarrow} \, X_2  
\nonumber \\
R_4 &:& \qquad\qquad X_2 \, \overset{K}{\longrightarrow} \, X_1.  
\nonumber 
\end{eqnarray}
We will study this system for large values of parameter $K \gg k_1V + k_2.$ 
Then reactions $R_3$
and $R_4$, occur many times on a timescale for which reactions $R_1$ and 
$R_2$ are unlikely to happen at all. In
such a regime, one might consider using multiscale methods to reduce
the computational cost of analysing the system. The slow quantity in
this system is $S = X_1 + X_2$. Note that this quantity is invariant
with respect to the fast reactions, and so only changes when either of
slow reactions ($R_1$ or $R_2$) occur.

The analytical solution of the steady state CME is given by the following 
multivariate Poisson distribution~\cite{Jahnke:2007:SCM}:
\begin{equation}\label{eq:poisson}
\mathbb{P} (X_1 = x_1, X_2 = x_2) =
\frac{\lambda_1^{x_1}\lambda_2^{x_2}}{x_1!x_2!} 
\, \exp \! \big[
- (\lambda_1 + \lambda_2)
\big],
\end{equation}
where $\lambda_1 = k_1 V / k_2$ and
$\lambda_2 = \lambda_1 (K+k_2)/ K.$ Let $\mathcal{P}(\lambda)$ be 
the Poisson distribution with rate $\lambda$ which is defined by
its probability mass function
\begin{equation*}
\mathbb{P} (X = x) = \frac{\lambda^{x}}{x!} \exp (- \lambda).
\end{equation*}
Using (\ref{eq:poisson}), we obtain that the exact
distribution of the slow variable $S = X_1 + X_2$ is 
\begin{equation}
\label{eq:slow}
\mathcal{P}(\lambda_0),
\quad
\mbox{where}
\quad
\lambda_0 = \lambda_1 + \lambda_2 = \frac{k_1 V}{k_2} 
\left ( 2 + \frac{k_2}{K} \right).
\end{equation}
In the rest of this section, we use (\ref{eq:slow}) to compare 
the accuracy and efficiency of the CMA, NMA and SSMA.
Each of the three algorithms gives us a different method to
approximate the effective drift and diffusion of the slow variable at
a given point on the slow manifold. For each method there are several
sources of error, and in this section we aim to identify the effect of
each, for each method.

\subsection{Quasi-steady state assumption error}\label{sec:qssa}
\begin{figure}[thb]
\centerline{
\raise 4.685cm \hbox{\raise 0.9mm \hbox{(a)}}
\hskip -5.6mm
\epsfig{file=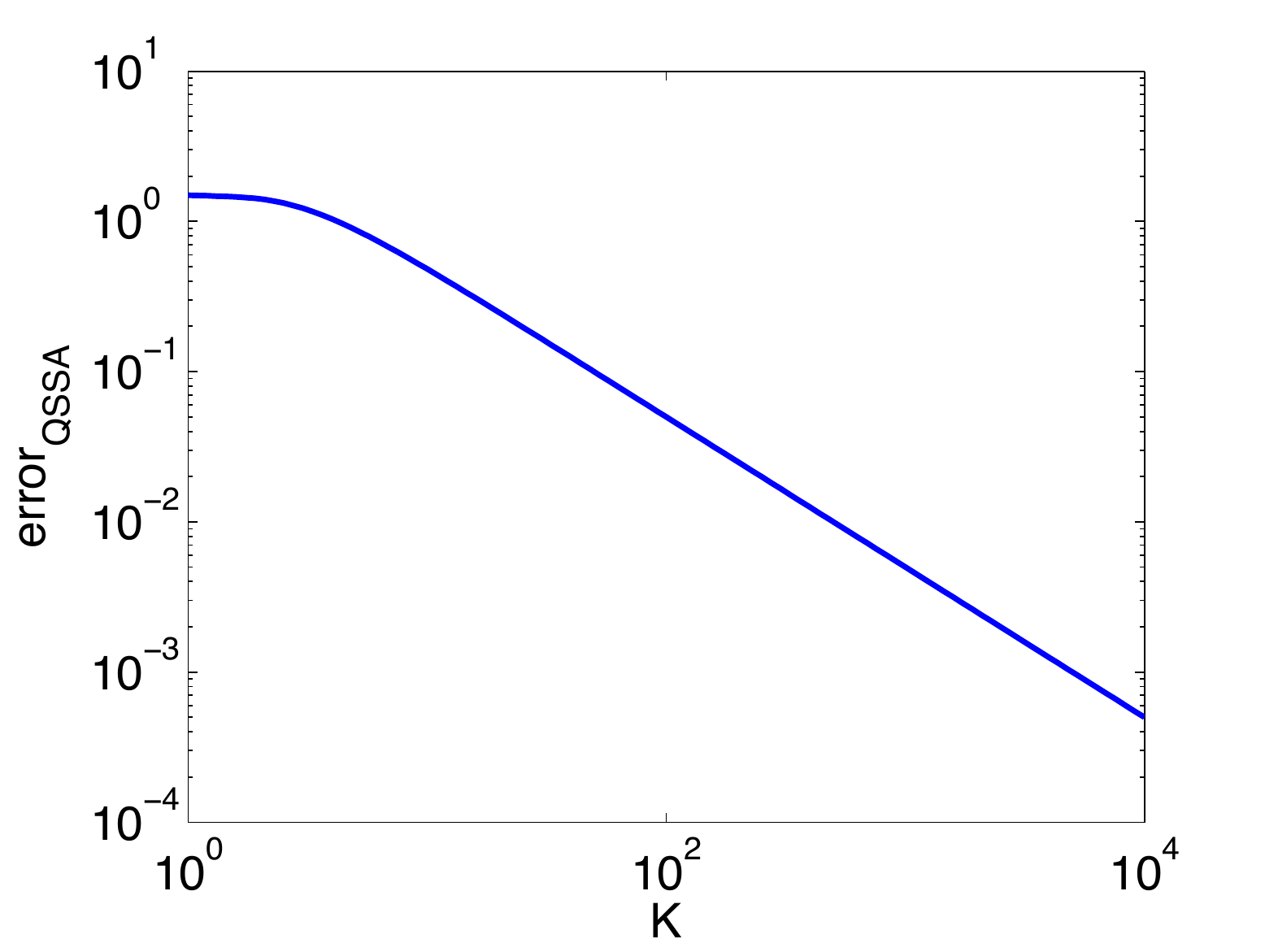,height=4.685cm}
\hskip 2mm
\raise 4.685cm \hbox{\raise 0.9mm \hbox{(b)}}
\hskip -5.6mm
\epsfig{file=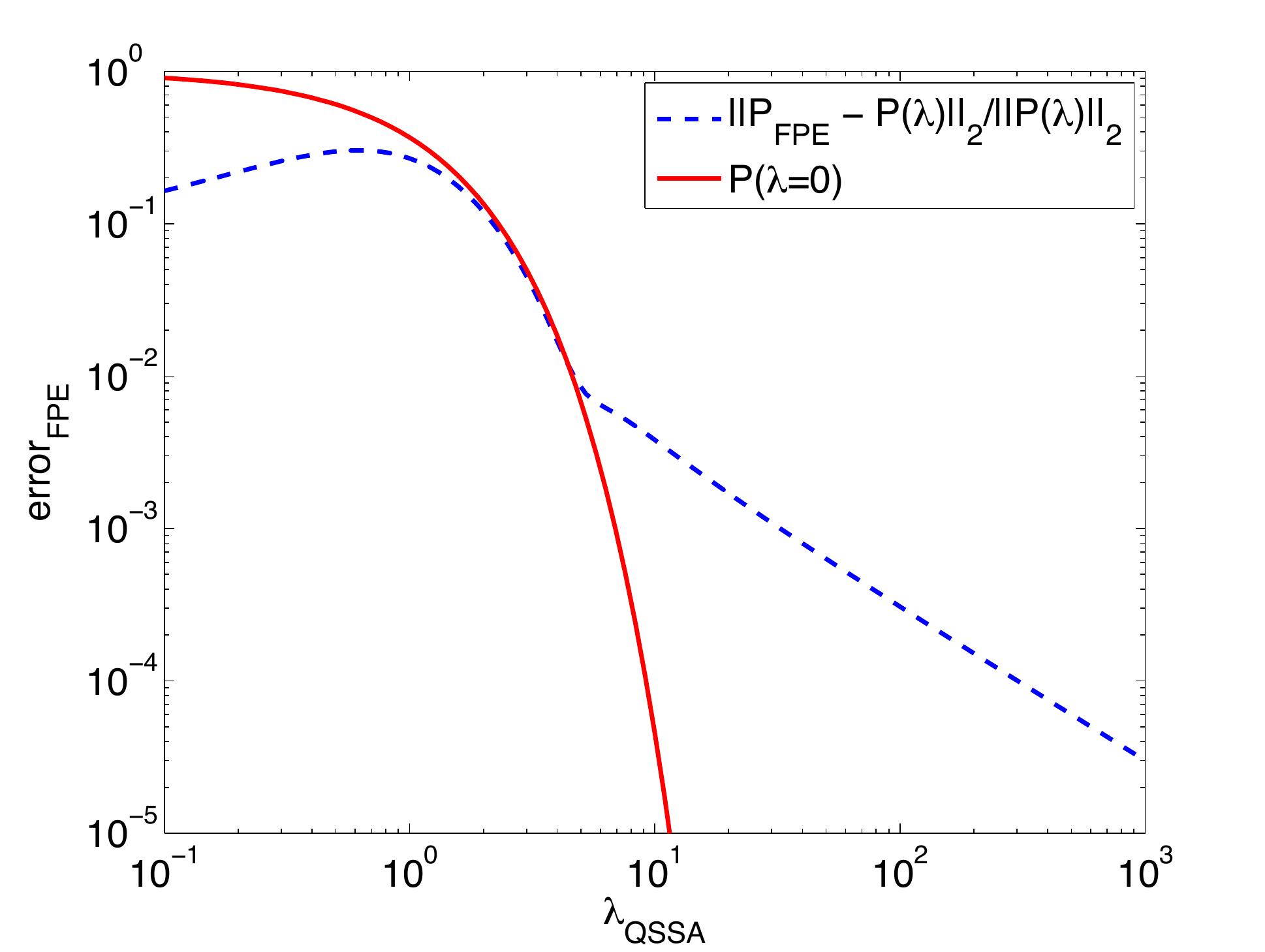,height=4.685cm}
}
\caption{{\rm (a)} {\it Plot of error $(\ref{eq:err})$, incurred by the QSSA 
  for the system $\eqref{eq:lin}$ with reaction parameters
  given $k_1 = k_2 = 1$ and $V=100$, as a function of $K$.} \hfill\break
  {\rm (b)} {\it Blue plot shows error $(\ref{errorFPE})$ as a function of
  $\qlam$, for the system $(\ref{eq:eff})$. Red plot shows value of
  $\mathcal{P}(\qlam)(0)$ as a function of $\qlam$.}
  \label{figure1}}
\end{figure}

The NMA and QSSMA both assume that the reactions can be partitioned into
fast and slow reactions. Both of these methods rely on the assumption
that the fast reactions enter equilibrium on a much faster (or even
instantaneous) timescale in comparison with the slow reactions. This
assumption leads to the approximation that the dynamics of the slow
variables can be described well by a system consisting only of the
slow reactions. For example, we assume that the variable $S = X_1 +
X_2$ in the system \eqref{eq:lin} can be well approximated by the
system:
\begin{equation}\label{eq:eff}
\LRARR{\emptyset}{S}{\bar{k}_1}{\bar{k}_2}.
\end{equation}
The two methods (NMA and QSSMA) differ in their calculation of the effective
reaction rates, $\bar{k}_1$ and $\bar{k}_2$. We denote the effective 
propensities for these two reactions $\bar{\alpha}_1(s)$ and 
$\bar{\alpha}_2(s)$ respectively. We will now isolate the error 
that is incurred by approximation of the full system by the reduced
system written in terms of the slow variables. Slow reaction $R_1$ does not depend on the 
value of the fast variables. Consequently, we have 
$$
\bar{\alpha}_1(s) = k_1 V
\qquad
\mbox{and}
\qquad
\bar{k}_1 = k_1. 
$$
The second effective rate $\bar{k}_2$ in (\ref{eq:eff}) has
to be calculated, because reaction $R_2$ includes fast
variables. The average values of $X_1$ and $X_2$ for the fast system 
(reactions $R_3$ and $R_4$ in \eqref{eq:lin}), for a given value of 
$S = s$, is~\cite{Cao:2005:MSS,Xue:2012:CMC}
$$
\mathbb{E} \big( X_1 | S = s \big) 
=
\mathbb{E} \big( X_2 | S = s \big)
=
\frac{s}{2}. 
$$
Therefore, we have 
$$
\bar{\alpha}_2(s) 
= 
k_2 \, \mathbb{E} \big( X_2 | S = s \big)
= \frac{k_2 \, s}{2} 
\qquad
\mbox{and}
\qquad
\bar{k}_2 
=
\frac{k_2}{2}.
$$ 
The probability density of $S$ is then given by the Poisson distribution
\begin{equation}
\label{eq:slowapprox}
\mathcal{P}(\qlam),
\quad
\mbox{where}
\quad
\qlam = \frac{2 \, k_1 V}{k_2}.
\end{equation}
We define the error incurred by the QSSA by
\begin{equation}\label{eq:err}
\mbox{error}_{\scriptscriptstyle \rm QSSA} 
=
\frac{\|\mathcal{P}(\qlam) 
- \mathcal{P}(\lambda_0)\|_2}{\|\mathcal{P}(\qlam)\|_2}.
\end{equation}
Comparing (\ref{eq:slow}) and (\ref{eq:slowapprox}),
we have
$$
\lambda_0 = \qlam + \frac{k_1 V}{K}.
$$
Therefore error (\ref{eq:err}) can be approximated
for large $K \gg 1$ by
$$
\mbox{error}_{\scriptscriptstyle \rm QSSA} 
=
\frac{1}{\|\mathcal{P}(\qlam)\|_2}
\sqrt{\sum_{n=0}^\infty
\left(
\frac{\qlam^{x}}{x!} \exp (- \qlam)
-
\frac{\lambda_0^{x}}{x!} \exp (- \lambda_0)
\right)^2}
\approx O(K^{-1}).
$$
Figure \ref{figure1}(a) shows 
error (\ref{eq:err}) a function of $K$. This plot confirms that 
this error decays like $K^{-1}$ as $K$ increases (gradient of 
linear part of plot is equal to $-0.998$).

The limitations of the stochastic quasi-steady-state approximation are
looked at in detail in~\cite{Thomas:2011:LSQ}.

\subsection{Diffusion Approximation Error}\label{sec:diff}
One of the sources of the error, common across all of three
methods (CMA, NMA, QSSMA), is that we are approximating a Markov 
jump process which has a discrete state space (the non-negative 
integers), by a SDE with a continuous state space (the positive
real numbers). To see the effect of this, let us consider the simple
birth-death chemical system (\ref{eq:eff}).
The steady state solution to the CME for this
system is given by the Poisson distribution (\ref{eq:slowapprox}).
The closest approximation that we can get to this process with 
an SDE, is the chemical Langevin equation~\cite{Gillespie:2000:CLE}. 
The corresponding stationary FPE for
this system is
\begin{equation}
\frac{1}{2} \frac{\partial^2}{\partial s^2}
\left( \Big( 
\bar{k}_1 + \bar{k}_2 s
\Big)
p(s)
\right)
-
\frac{\partial}{\partial s}
\left( \Big( \bar{k}_1 - \bar{k}_2 s \Big)
p(s)
\right)
=
0.
\label{FPEBD}
\end{equation}
It can be explicitly solved to get~\cite{Sjoberg:2009:FPA}.
\begin{equation}\label{eq:density}
p(s) =
  C\exp(-2s)(s + \qlam)^{4 \qlam - 1},
\end{equation}
where $C$ is determined by normalisation 
$\int p(s) \, \mbox{d} s = 1.$
Here the integral is assumed to be taken over $s \ge 0$ (if we 
want to interpret $s$ as concentration) or $s \ge - \qlam$
(if we do not want to impose artificial boundary conditions
at $s=0$). 
Considering more complicated systems, it is more natural 
to assume that the domain of the chemical Langevin equation is
a complex plane~\cite{Schnoerr:2014:CCL}.

Suppose we now wish to quantify the difference between
probability distributions (\ref{eq:slowapprox})
and (\ref{eq:density}) as a function of $\qlam$. The first 
issue that we come across is that the solution of the steady-state CME 
is a probability mass function, and the solution to the steady-state
FPE is a probability density function.
However, we can simply integrate the probability density function 
over an interval of size one centered on each integer, to project 
this distribution onto a discrete state space with mass function 
$\pFPE$ so that the two distributions can be compared:
$$
\pFPE(n) = \int^{n+1/2}_{n-1/2} P(s) \, \mbox{d}s. 
$$
Another issue is to specify the treatment of negative values of $s$.
In our case we truncate the distribution for $s \ge 0$.
We can then exactly calculate $\pFPE$ to get
$$
\pFPE(n) 
=
\frac{1}{\Gamma(4\qlam) - \gamma(4\qlam,2\qlam)}
\Big( 
\gamma(4\qlam, 2 n + 1 + 2 \qlam)
-
\gamma(4\qlam, 2 n - 1 + 2 \qlam)
\Big),
$$
where $\gamma(k,x) = \int_0^x z^{k-1}\exp(-z) \,\mbox{d}z$ denotes the
lower incomplete gamma function and 
$\Gamma(k) = \int_0^\infty z^{k-1}\exp(-z) \, \mbox{d} z$ is the 
gamma function. Then we can consider the $l^2$ difference between these two
distributions for a given value of $\qlam$,
\begin{equation}
\mbox{error}_{FPE}
=
\frac{\|\pFPE - \mathcal{P}(\qlam)\|_2 }{\| \mathcal{P}(\qlam)\|_2}
= 
\sqrt{
\sum_{n=0}^\infty \left(
\pFPE(n) 
- 
\frac{\qlam^n}{n!} \, \exp[ - \qlam ] \right)^2
}.
\label{errorFPE}
\end{equation}
Figure \ref{figure1}(b) shows how $\mbox{error}_{FPE}$
decays as $\qlam$ increases. The slightly odd sickle
shaped error curve for small $\qlam$ is due to the probability mass
of $\mathcal{P}(\qlam)$ being peaked close to (or at) zero. 
In this region the
diffusion approximation is very poor. To illustrate this, the value of
$\mathcal{P}(\qlam)$ at $s=0$, i.e.
$\exp[ - \qlam ]$,
is also plotted in Figure \ref{figure1}(b) (red
curve). Once the peak of the probability mass has moved far enough
away from zero, then the error is no longer dominated by being too
close to zero, and decays inversely proportional to
$\qlam$ (gradient of log-log plot is $-1.044$).

In the 2000 paper by Gillespie~\cite{Gillespie:2000:CLE} on the
chemical Langevin equation, a condition is put on the types of system
which are well approximated by a diffusion. The probability of the
system entering a state where the copy numbers of one or more of the
chemical species in the state vector are close to zero must be small. Otherwise the
approximation becomes poor.
In the case of diffusion approximations of the slow variable(s) of a
system, the trajectories must likewise stay away from regions of the
state-space with low values of the slow variable(s).
Further discussions and results regarding the accuracy of the chemical
Langevin and Fokker-Planck equations can be found in~\cite{Grima:2011:ACL}.

\subsection{Monte Carlo Error}\label{sec:MC}
The CMA and NMA both employ bursts of stochastic simulations
to estimate the effective drift and diffusion of the slow variable. 
The main advantage of the QSSMA is that no such computation is required.

In the case of the CMA, as described in Section \ref{sec:CMA}, the
full system is simulated, including fast and slow reactions. The
computed trajectory is constrained to a particular value of the slow
variables, and so whenever a slow reaction occurs, the trajectory is
projected back onto the same point on the slow manifold. The effective
drift and diffusion at that point on the slow manifold are functions
of statistics about the relative frequency of these slow reactions. As
such, as given by the central limit theorem, the error in these
estimations are mean zero and normally distributed, with variance
proportional to $N_{S}^{-1}$, where $N_{S}$ is the number of slow
reactions simulated during the estimation. Since it is necessary to
simulate all of the fast and slow reactions in the system, depending
on the difference in timescales this can be very costly. Since 
the ratio of occurrences of fast reactions to slow reactions is
proportional to $K$, the cost of the estimation is
$\mathcal{O}(K \, N_{S})$.

In comparison, the NMA, as described in Section \ref{sec:NMA}, is
only concerned with finding the average value of the fast variables
through a Gillespie SSA simulation of only the fast
variables. Therefore, the Monte Carlo error is again mean zero and
normally distributed, with variance $N_{F}^{-1}$, where $N_{F}$ is the
number of fast reactions simulated. Since we only simulate the fast
reactions, the cost of the estimation is
$\mathcal{O}(N_{F})$. Therefore, the cost of estimation for the CMA is
approximately $K$-times that of the NMA for the same degree
of accuracy.

\subsection{Approximation of the solution of the stationary FPE}
All three of the algorithms (CMA, NMA, QSSMA) also incur error 
through discretisation errors in the approximation of the solution 
to the steady-state FPE \eqref{eq:FPE}. The error of such an
approximation is dependent on the method used, such as the adaptive
finite element method used for three-dimensional chemical systems in~\cite{Cotter:2013:AFE}. In this paper we are mainly interested in 
the accuracy of the various methods for estimating the
effective drift and diffusion functions, and as such we aim to
simplify the methods for solving this PDE as much as
possible. Therefore we will only consider systems with one-dimensional
slow variables. The steady state equation corresponding to 
FPE \eqref{eq:FPE} is then effectively an ordinary differential 
equation, and one which can be solved directly~\cite{Gillespie:2000:CLE}
to obtain
\begin{equation}
\label{FPSol}
\mathcal{P}_{S}(s) = \frac{C}{D(s)} \exp 
\left( 
  \int_0^s \frac{V(z)}{D(z)} \, \mbox{d}z \right),
\end{equation} 
where $C$ is the normalisation constant.
With approximations of $V$ and $D$ over a range of values of $s$
(through implementation of the CMA, NMA or QSSMA), the integral in
this expression can be evaluated using a standard quadrature routine
(for instance the trapezoidal rule). The errors incurred here will 
be proportional to the grid size to some power, depending on the 
approximation method used.

\subsection{Comparison of Sources of Error}

\begin{table}[t]
\centering
\begin{tabular}{|c|c|c|c|}
\hline
Source of Error & CMA & NMA & QSSMA \\ \hline
Diffusion approximation & \tick & \tick & \tick \\ 
QSSA & \cross & \tick & \tick \\ 
Monte Carlo error & Cost - $\mathcal{O}(\mathcal{K}N_S)$ & 
Cost - $\mathcal{O}(N_F)$ & \cross \\
PDE approximation & \tick & \tick & \tick \\
\hline
\end{tabular}
\caption{\it Table to summarise the sources of error of the 
three algorithms.\label{tab:errors}}
\end{table}

Table \ref{tab:errors} summarises the analysis of the errors of the
various methods that we have looked at in this section. Each method
has advantages and disadvantages, depending on the type of system
which a modeller wishes to apply the methods to. All of the methods use
diffusion approximations, and as such, if the slow variables of the 
system of interest cannot be well approximated by a diffusion, then 
none of the proposed methods are suitable. If the QSSA does not
hold for the system, then the CMA is the best choice. If
it does hold, and the analytical solution for the effective
propensities is available to us, then the QSSMA is the best choice,
since it does not incur Monte Carlo, and is the least expensive of the
three algorithms. Finally, if no such analytical solution is
available, but the QSSA still holds, then the
NMA is the best choice of algorithm, since it converges faster than
the CMA.
 
Next, we apply the three methods to three different
parameter regimes of the system given by \eqref{eq:lin}. In each of
the experiments, we set $k_1 = k_2 = 1$ and $V=100$
and we vary $K$. We use $K=10$, $K=200$ and $K=10^3$.
In each case, the CMA, NMA and the
QSSMA are all applied to the system over the range of values $S = X_1
+ X_2 \in [101,300]$. This range is chosen since the
invariant distribution of the slow variable (\ref{eq:slow})
is the Poisson random variable with intensity 
$\lambda_0 = 200 + 100/K$, and therefore the vast majority of the invariant
density is contained with this region for all three parameter
regimes. Furthermore, we implemented these algorithms on a computer with 
four cores, and to optimise the efficiency of parallelisation it was
simplest to choose a domain with the number of states divisible by 4,
hence the region starting at 101 as opposed to 100. For the CMA and
the NMA, a range of different numbers of reactions were used in order
to estimate the drift and diffusion parameters at each point, $N_S,
N_F \in \{10^1, 10^2,\ldots,10^9,10^{10}\}$.
Each code was implemented in C, and optimised to the best of the
authors' ability, although faster codes are bound to be possible. The
number of CPU cycles used was counted using the C time library. The
CPU cycles used over all 4 cores were collated into a single number
for each experiment.

\begin{figure}[t]
\centerline{
\raise 4.685cm \hbox{\raise 0.9mm \hbox{(a)}}
\hskip -5.6mm
\epsfig{file=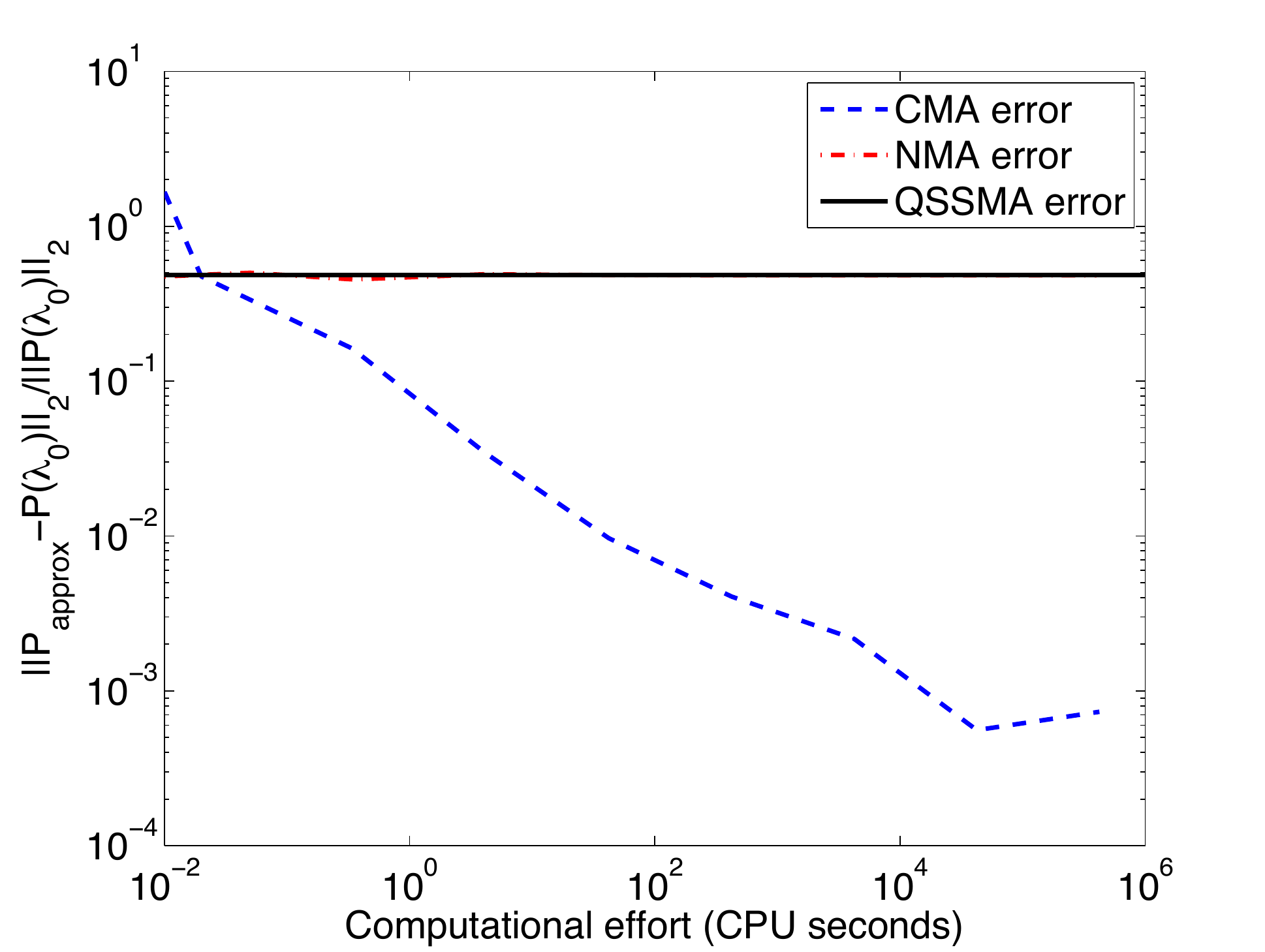,height=4.685cm}
\hskip 2mm
\raise 4.685cm \hbox{\raise 0.9mm \hbox{(b)}}
\hskip -5.6mm
\epsfig{file=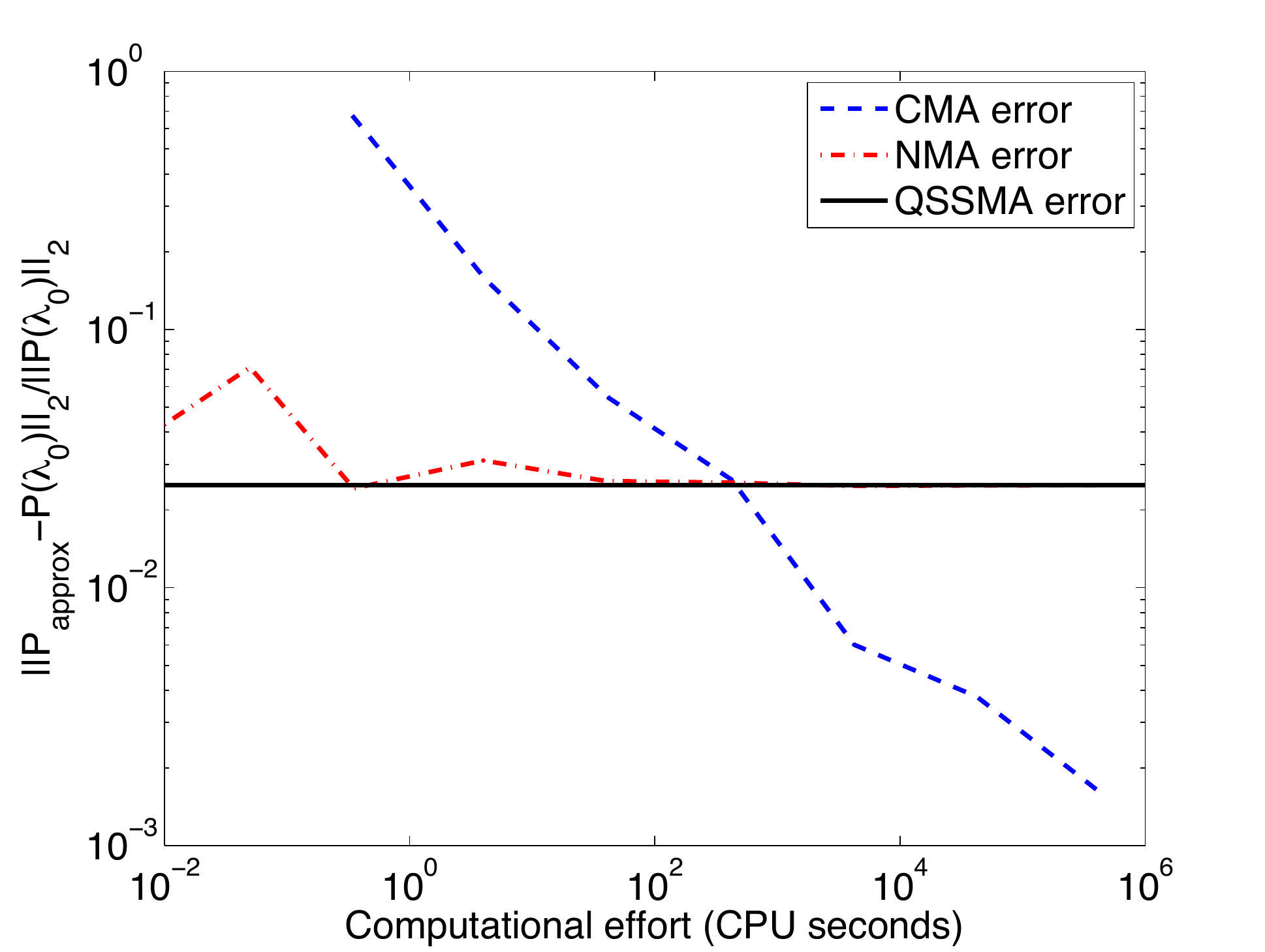,height=4.685cm}
} 
\vskip 2mm
\leftline{
\raise 4.685cm \hbox{\raise 0.9mm \hbox{(c)}}
\hskip -5.6mm
\epsfig{file=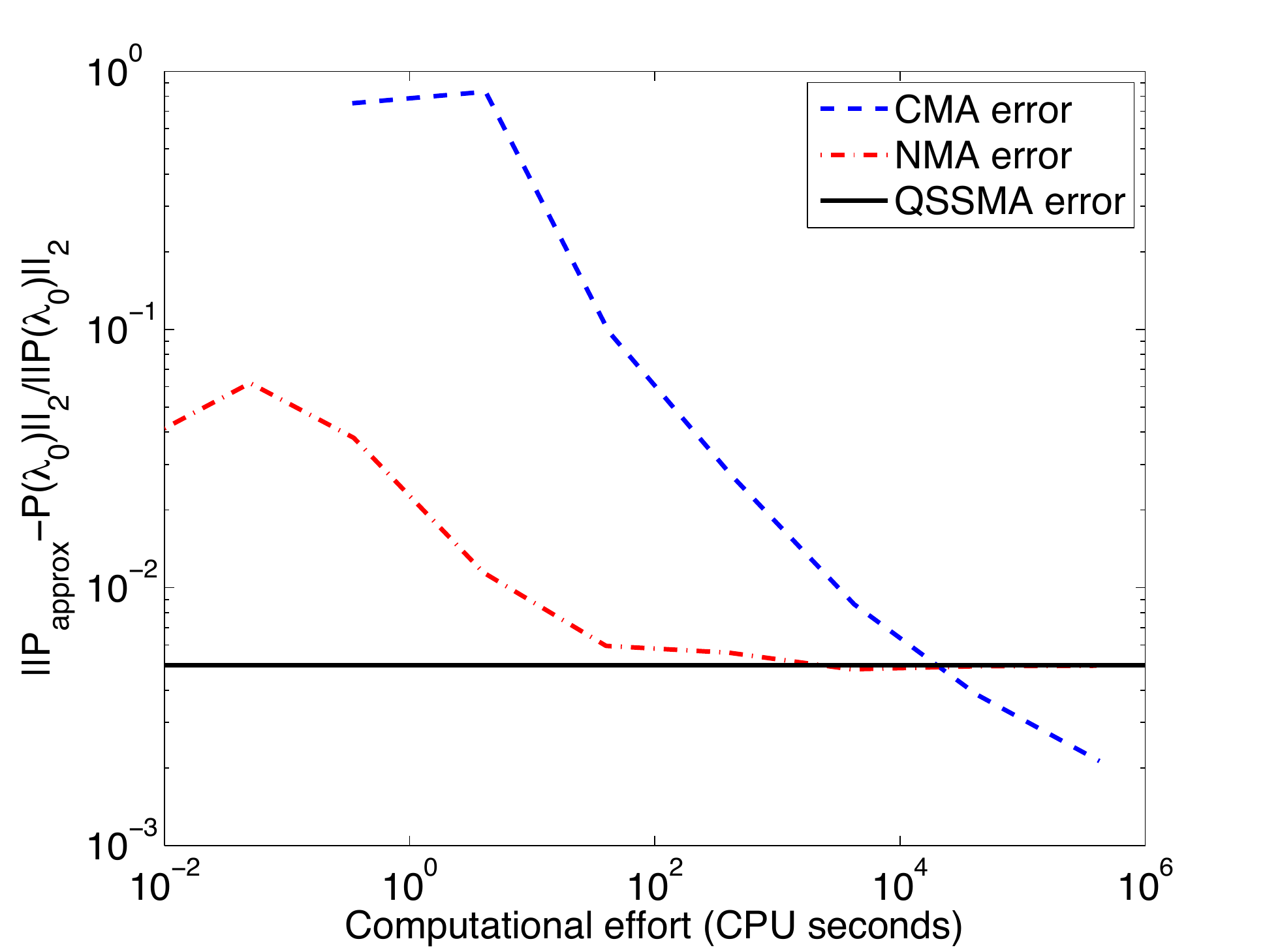,height=4.685cm}
} 
\caption{{\rm (a)}
  {\it Error $(\ref{totalerror})$ of the CMA (blue), NMA (red) and the QSSMA
  (black) as the function of the computational effort   
  for the chemical system $(\ref{eq:lin})$ with parameter values
  $k_1 = k_2 = 1$, $V=100$ and $K = 10$. For illustrative
  purposes the single value of the QSSMA error is plotted as a
  horizontal line. Actual computation time for the QSSMA is
  negligible.} \hfill\break
  {\rm (b)} {\it As in panel} {\rm (a)}, {\it with $K = 200.$} \hfill\break 
  {\rm (c)} {\it As in panel} {\rm (a)}, {\it with $K = 10^3.$} 
  \hfill\break.
  \label{figure2}}
\end{figure}

The results of these experiments are shown in Figure \ref{figure2}.
Note that the results of the QSSMA and NMA are unaffected by a change
in the value of $K$. In the case of the NMA, a change in the
value of this variable simply scales time in computation of the fast
subsystem, but does not affect the result. As such, only one run was
necessary for these methods for the three different parameter
regimes. The three plots for the NMA use the same simulations, but
since the true distribution of the slow variable is affected by the
change in $K$, the error plots differ across the three parameter
regimes. In Figure \ref{figure2} we plot the total relative error
\begin{equation}
\mbox{error}
=
\frac{\| P_{approx}
- \mathcal{P}(\lambda_0)\|_2}{\|\mathcal{P}(\lambda_0)\|_2},
\label{totalerror}
\end{equation}
where $P_{approx}$ is the steady state distribution of the slow
variable obtained by the corresponding method and
$\mathcal{P}(\lambda_0)$ is the exact solution given by
(\ref{eq:slow}). Figure \ref{figure2}(a), with $K = 10$,
denotes a parameter regime in which the QSSA produces a
great deal of error. The NMA error quickly converges to the level seen
in the QSSMA, but neither method can improve beyond this. The error
seen from the CMA improves on both of these methods with relatively
little computational effort. One might argue that the system is not
actually a ``true'' multiscale system in this parameter regime, but
the CMA still represents a good method for analysis of the dynamics of
such a system, since its implementation can be parallelised in a way
which scales linearly with the number of cores used.

Figure \ref{figure2}(c) shows a parameter regime which is highly
multiscale. In this regime, the QSSA is far more reasonable,
and as such we see vastly better performance from the NMA and QSSMA
methods. However, eventually the CMA still has a higher accuracy 
than these two other
approaches, although at not inconsiderable computational cost. In this
case, the error incurred by the QSSA may be considered small
enough that a modeller may be satisfied enough to use the QSSMA, whose
costs are negligible. For more complex systems, the CME
for the fast subsystem may not be exactly solvable, or even easily
approximated, and in these cases the NMA would be an appropriate 
choice. If a more accurate approximation is required, once again 
the CMA could be used. In summary, even for simple system (\ref{eq:lin}), 
with different parameter regimes, different methods could be considered 
to be preferable. 

\section{A bistable example}
\label{bistsection}
\begin{figure}[t]
\centering
\includegraphics[width=0.6\textwidth]{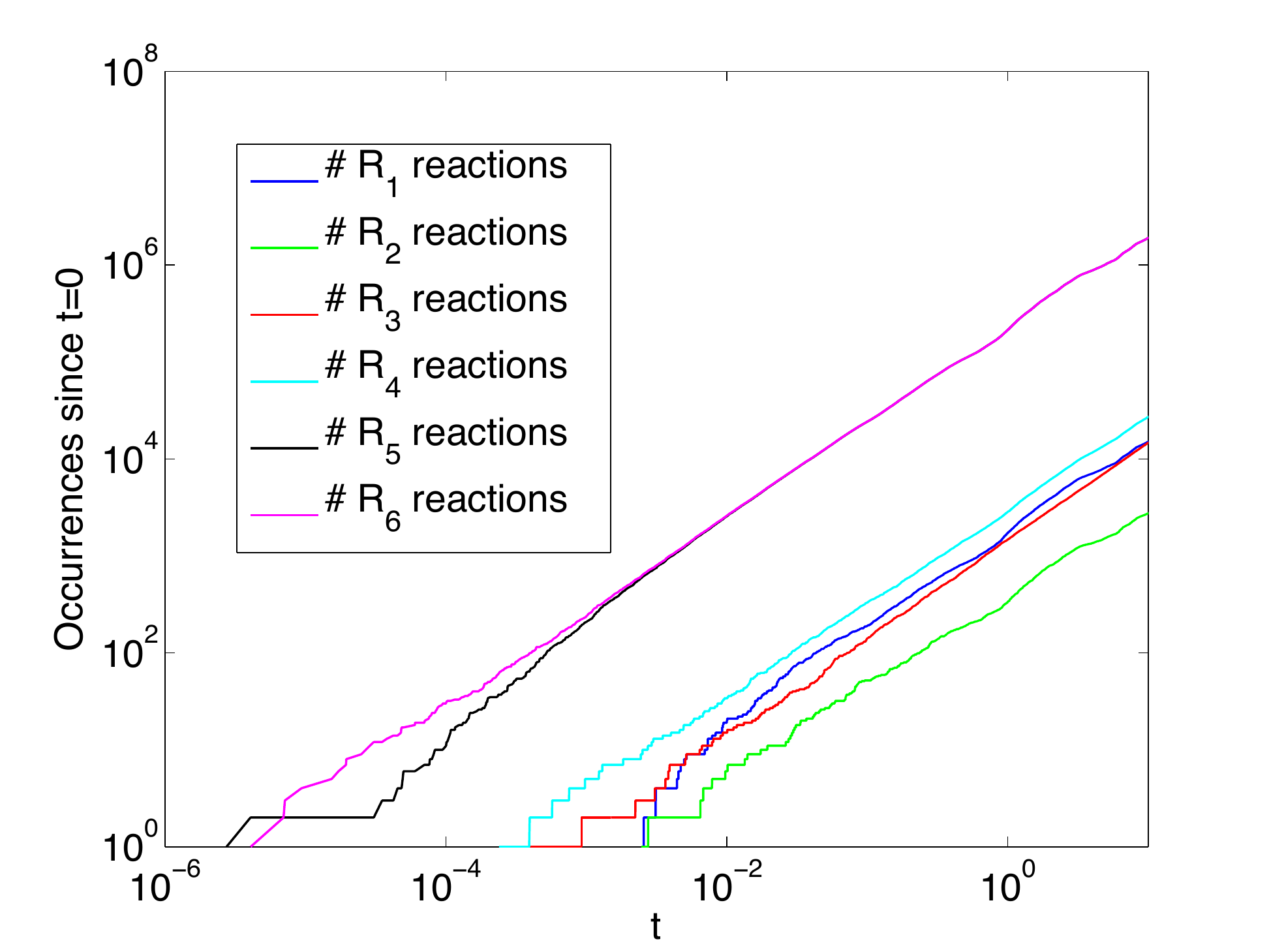}
\caption{{\it Plot to show the frequency of the different 
reactions in the
  system \eqref{eq:bis} with parameters given by 
  \eqref{eq:params:bis}, computed using
  the Gillespie SSA, with initial condition ${\mathbf X} = [100,100]$.}
  \label{figure3}}
\end{figure}

In this section, we will compare the presented methods for a
multiscale chemical system which is bimodal, 
with trajectories switching between
two highly favourable regions~\cite{Cotter:2011:CAM}:
\begin{equation}
\LRARR{X_2}{X_1+X_2}{k_{1}}{k_{2}},\qquad 
\LRARR{\emptyset}{X_1}{k_{3}}{k_{4}}, \qquad 
\LRARR{X_1+X_1}{X_2}{k_{5}}{k_{6}}.
\label{eq:bis}
\end{equation}
One example of the parameter values for which this system is bistable
for is given by:
\begin{equation}
k_1 = 32, \quad \frac{k_2}{V} = 0.04, \quad k_3V  = 1475,\quad
k_4 = 19.75, \quad \frac{k_5}{V}= 10, \quad k_6 = 4000.
\label{eq:params:bis}
\end{equation}
In this regime, reactions $R_5$ and $R_6$ (with rates $k_5$ and $k_6$
respectively) are occurring on a much faster time scale than the
others. This can be seen in Figure \ref{figure3}, which shows
the cumulative number of occurrences of each reaction in this system,
simulated using the Gillespie SSA, given in Table \ref{SSAtable}.
Both the species $X_1$ and $X_2$ are fast variables, since neither is
invariant to the fast reactions. As in~\cite{Cotter:2011:CAM}, we pick a slow variable which is
invariant to reactions $R_5$ and $R_6$, $S = X_1 + 2X_2$. We now aim
to compare the efficiency and accuracy of the CMA, NMA and QSSMA
in approximating the stationary distribution in two different
parameter regimes, with different spectral gaps between the fast and
slow variables. This can be done by altering the values of rates $k_5$
and $k_6$.
 
One issue with using such a system, is that we cannot compute the
analytical expression for the invariant density. Therefore, we compare
with an accurate approximation of the invariant density.

\subsection{Application of the NMA, QSSMA and CMA to reaction system \eqref{eq:bis}}\label{sec:new}
For the NMA and the QSSMA, we assume that reactions $R_5$ and $R_6$
are occurring on a very fast time scale. Therefore we may assume that between
reactions of other types (the slow reactions $R_1$-$R_4$), the
subsystem involving only the fast reactions enters statistical
equilibrium. In this case, we can reformulate \eqref{eq:bis} such that
it is effectively a set of reactions which change the slow variable
(reactions $R_5$ and $R_6$ are of course omitted as the slow 
variable $S=X_1+2X_2$
is invariant with respect to these reactions):
\begin{equation}\label{eq:QSSbis}
\LRARR{2S}{3S}{\widehat{k}_{1}}{\widehat{k}_{2}},\qquad 
\mbox{\raise 0.4mm \hbox{and}} \qquad
\LRARR{\emptyset}{S}{\widehat{k}_{3}}{\widehat{k}_{4}}.
\end{equation}
Next, we design the QSSMA by analytically computing the rates 
$\widehat{k}_{1}$, $\widehat{k}_{2}$ , $\widehat{k}_{3}$
and $\widehat{k}_{4}$.

To estimate the rates in (\ref{eq:QSSbis}) 
we must compute the average values of $X_1$ and $X_2$ 
for fixed value of $S=s$ in the fast
subsystem:
\begin{equation}\label{eq:FSS}
\LRARR{X_1+X_1}{X_2}{k_5}{k_6},\qquad X_1(0) + 2X_2(0) = s.
\end{equation}
As in~\cite{Cao:2005:SSS}, we approximate the mean values 
$\langle X_1 \rangle$
and $\langle X_2 \rangle$ of $X_1$ and $X_2$, respectively, by the solutions
of the deterministic reaction rates equations. The authors showed
that such an approximation is reasonably accurate for this particular
fast subsystem (a reversible dimerisation reaction). Therefore, we
need to solve the following system:
\begin{eqnarray*}
\frac{k_5}{V} \langle X_1 \rangle^2 
= k_6 \langle X_2 \rangle, \qquad  
\langle X_1 \rangle + 2 \langle X_2 \rangle = s. 
\end{eqnarray*}
The nonnegative unique solution is given by:
\begin{equation}\label{eq:RRES}
\langle X_1 \rangle = \frac{Vk_6}{4k_5} \left ( \sqrt{1 + 8\frac{k_5}{Vk_6}s} - 1
\right)
\qquad \mbox{and} \qquad
\langle X_2 \rangle =   \frac{s - \langle X_1 \rangle }{2}.
\end{equation}
Then the effective propensity function of the slow reaction can be 
written as~\cite{Cao:2005:SSS}:
\begin{eqnarray}
\widehat{\alpha}_1 &=& k_1 \, \langle X_2 \rangle,\nonumber\\
\widehat{\alpha}_2 &=& 
\frac{k_2 \, s \, \langle X_2 \rangle}{V}  
- 
\frac{2 \, k_2 \, \langle X_2 \rangle^2}{V}
+ 
\frac{2 \, k_2 \, k_6 \, \langle X_2 \rangle}{8 \, k_5
\langle X_2 \rangle - 2 \, k_5 (2s
  + 3) - k_6 \, V} \; ,
\label{eq:EP} \\
\nonumber
\widehat{\alpha}_3 &=& k_3 \, V, \\ 
\nonumber
\widehat{\alpha}_4 &=& k_4 \langle X_1 \rangle = k_4 \, \left( s - 2 \langle X_2 \rangle \right).
\end{eqnarray}
More computational effort is required for this problem, since we need
to compute $\langle X_2 \rangle$ for each value of $S=s$ on the mesh. 
However, the computational effort is still negligible in comparison 
with the CMA or NMA.
The biggest drawback with the QSSMA is the increase in mathematical
complexity as the fast and slow systems themselves become more
complicated. The more complexity there is, the more approximations
need to be made in order to find the values of the effective
propensities.
The NMA simulations of the fast subsystem \eqref{eq:FSS} in order to 
approximate the effective propensities \eqref{eq:EP}, which are
then fed into the Fokker-Planck equation for the slow sub-system.
For the CMA simulations, we let $S = X_1 + 2X_2$ be the slow
variable, and we let $X_2$ be the fast variable. For further details
on how to apply the CMA to this system, see~\cite{Cotter:2011:CAM}.

\subsection{Numerical Results}
In general, systems which have second (or through modelling
assumptions, higher) order reactions cannot be solved exactly,
although there are some special cases which can be
tackled~\cite{Grima:2012:SSF}. System (\ref{eq:bis}) has second order reactions
and hence we assume that the invariant distribution for the system cannot be solved
analytically. 
As such, we are not able to compare the approximations
arising from the three methods considered (CMA, QSSMA, NMA) to an exact
solution. However, we can approximate the solution to the CME
for this system, as we did in~\cite{Cotter:2011:CAM},
by solving it on a truncated domain, assuming that the invariant
distribution has zero probability mass everywhere outside this
domain.

\begin{figure}[t]
\centerline{
\raise 4.685cm \hbox{\raise 0.9mm \hbox{(a)}}
\hskip -5.6mm
\epsfig{file=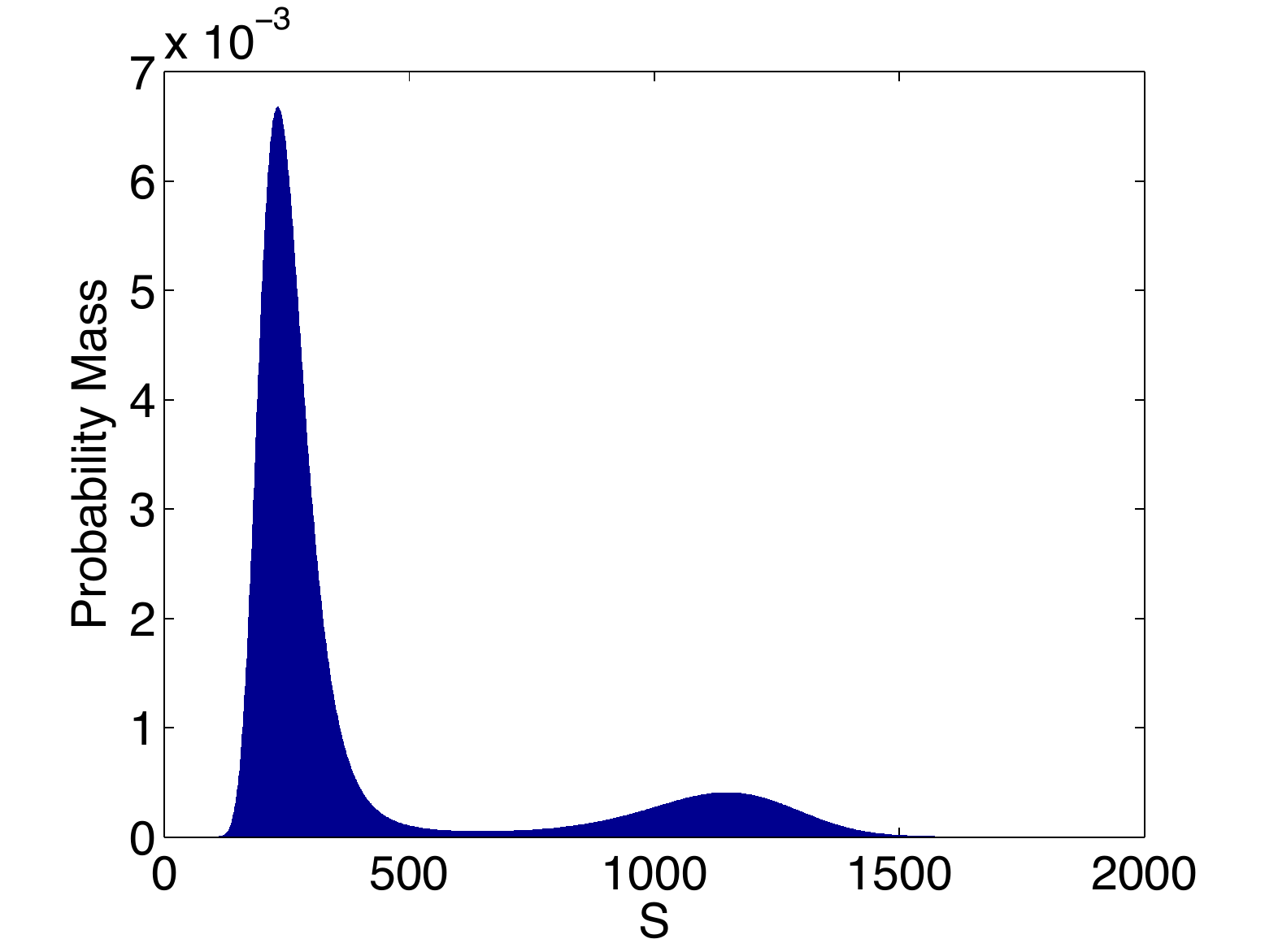,height=4.685cm}
\hskip 2mm
\raise 4.685cm \hbox{\raise 0.9mm \hbox{(b)}}
\hskip -5.6mm
\epsfig{file=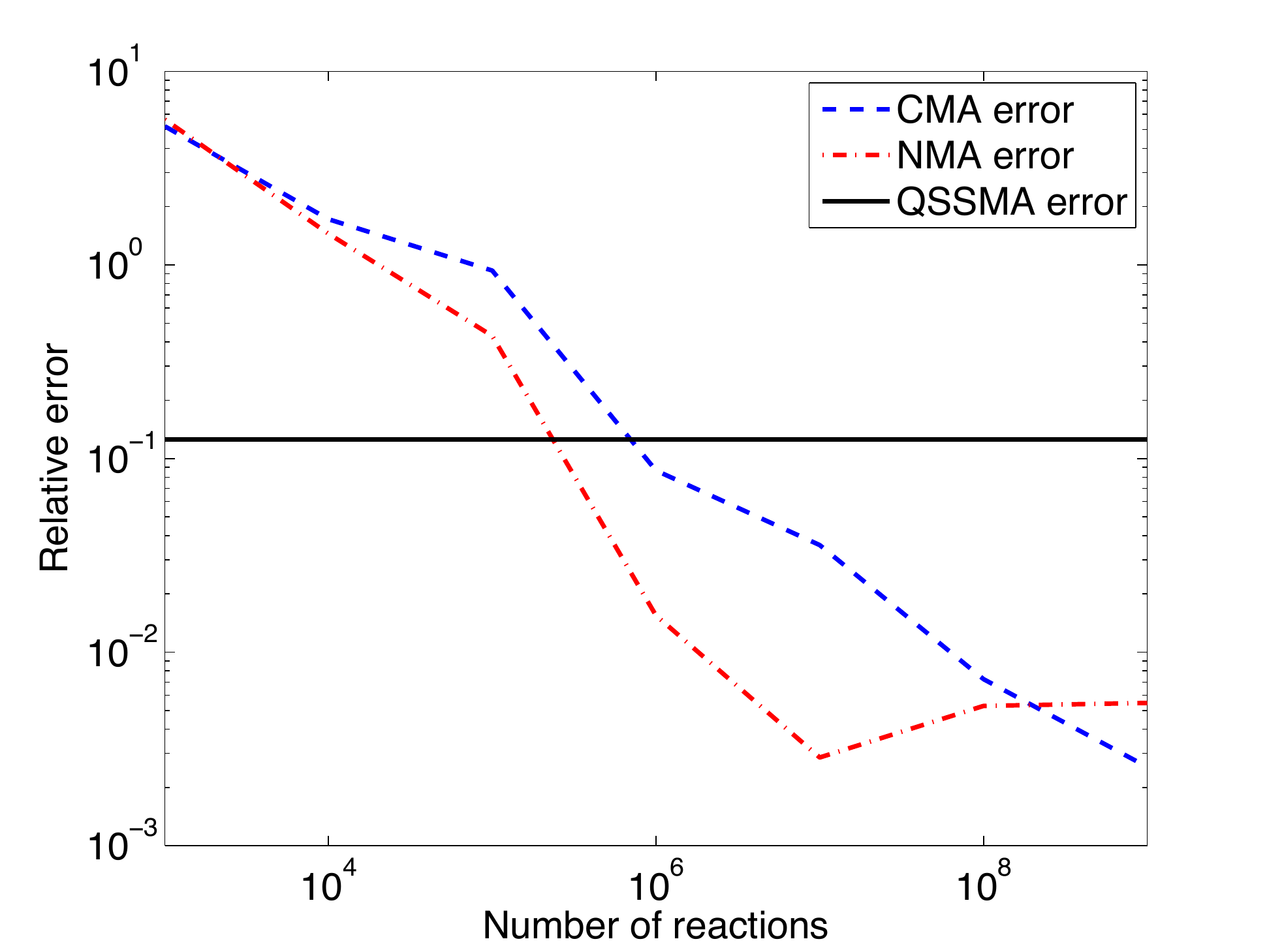,height=4.685cm}
} 
\caption{
{\rm (a)} 
{\it Approximation of the solution of the CME for the system
\eqref{eq:bis} with parameters \eqref{eq:params:bis}. 
The domain was truncated to $\Omega = [0,10^3] \times [0,1.5 \times 10^3]$.}
\hfill \break
{\rm (b)} {\it Approximation of the errors of the CMA, NMA and QSSMA 
for the system
  \eqref{eq:bis} with parameters \eqref{eq:params:bis}. Error was 
  estimated using the
  approximation of the solution of the CME shown in
  panel} {\rm (a).}
\label{figure4}}
\end{figure}

For the numerics that follow, we solve the CME on a truncated domain
$\Omega = [0,10^3] \times [0,1.5 \times 10^3]$. 
The CME on this domain gives us a
sparse linear system, whose null space is exactly one dimensional. The corresponding
null vector gives us, up to a constant of proportionality, our
approximation of the solution of the CME. We normalise this vector, and then sum the invariant probability mass over states
with the same value of the slow variable $S$. This procedure gives us
the approximation of the invariant density of the slow variable which
is plotted in Figure \ref{figure4}(a). Although this is only an
approximation, it is a very accurate one. To demonstrate this, we
compared the approximation of the solution $\mathcal{P}_\Omega$ of  the steady-state CME on our chosen domain
$\Omega = [0,10^3] \times [0,1.5 \times 10^3]$ with an approximation $\mathcal{P}_\omega$
over a smaller domain, $\omega = [0,8.0\times 10^2] \times [0,1.25 \times
10^3]$. The relative $l^2$ difference between these two approximations
was $1.4571\times10^{-11}$, indicating that any error in the
approximation $\mathcal{P}_\Omega$ is highly unlikely to be of the order of magnitude of
the errors seen in Figure \ref{figure4}(b), where we have used
$\mathcal{P}_\Omega$ to
approximate the error of multiscale methods.

Figure \ref{figure4}(b) shows the error plots for the three methods
for the system \eqref{eq:bis}, using the approximation of the solution of the CME as
the 
``exact'' solution. The computational effort here is measured in terms
of the numbers of simulated reactions. Since the computational
cost for the two methods which use Monte Carlo simulations are
dominated by the cost of producing two pseudo-random numbers for each
iteration, this is a good estimate. 

Unlike in the previous example, the distribution of the fast variables
for a fixed value of the slow variables is not known analytically, and
as such an approximation has been made, as was outlined in Section
\ref{sec:new}. As such, we no longer expect the NMA
approximation to converge to the QSSMA approximation as in the last
example. This can be seen in the error plot in Figure
\ref{figure4}(b), where the error in the QSSMA approximation, which
again had negligible cost to compute, is relatively high, with the NMA
and CMA quickly outperforming it. The NMA slightly
outperform the CMA at first, but appears to be unable to improve past
a relative error of around $7\times 10^{-2}$. Note that this is 9
orders of magnitude bigger than the relative $l^2$ difference between
$\mathcal{P}_\Omega$ and $\mathcal{P}_\omega$, and so
is highly unlikely to be an artefact of the method we have used to
approximate this error. As in the previous
example, and as seen in~\cite{Cotter:2011:CAM}, although the cost of
the CMA is in general higher than the other methods, if a high
precision solution is required, it is arguably the method of
choice, as the error continues to decrease monotonically 
$\sim \mathcal{O}(\sqrt{N_s})$.

\section{Discussion and conclusions}\label{sec:conc}

In this paper we have introduced two new methods for approximating the
dynamics of slow variables in multiscale stochastic chemical networks,
the NMA and QSSMA. These new methods combine ideas from the
CMA~\cite{Cotter:2011:CAM}, with ideas used in existing methods for
speeding up the Gillespie SSA for multiscale
systems~\cite{E:2005:NSS,Cao:2005:SSS}. We then undertook a detailed
numerical study of the different sources of error that these methods
incur, for a simple chemical system for which we have an exact
solution of the CME. Error is incurred due to the
approximation of the dynamics by a diffusion process, Monte Carlo
error in the approximation of the effective drift and diffusion terms,
error due to application of the QSSA, and
numerical error in approximation of the steady-state Fokker Planck
equation, in various combinations. We then conducted another numerical
study for a bistable chemical system, approximating the error by using
an approximation of the solution to the CME for
the system.

What we may surmise from this work, is that different methods,
utilising different types of approximations, are appropriate for
different types of system, or even in different parts of the parameter
space of the same system.  The methods in this paper differ from many
others for stochastic fast/slow systems mainly in the
approach of approximating the slow variables by a SDE. 
The majority of the other methods proposed in
the literature use different ways of speeding up the Gillespie SSA
for multiscale
systems~\cite{Rao:2003:SCK,Haseltine:2002:ASC,E:2005:NSS,Cao:2005:MSS}. In
each, the full simulation of the fast species is replaced by some sort
of approximation, so that an SSA-type algorithm, just for the slow
species, may be implemented.

Many other approaches in the literature rely on a QSSA, 
including that taken by Rao and Arkin~\cite{Rao:2003:SCK}, Haseltine and
Rawlings~\cite{Haseltine:2002:ASC}, E, Liu and
Vanden-Eijnden~\cite{E:2005:NSS} and Cao, Gillespie and
Petzold~\cite{Cao:2005:MSS}. All of these methods will incur the error
that is seen in Figure \ref{figure1}(a). This error will be
negligible for some systems, and significant for others, as we saw in
Section \ref{sec:qssa}, and in the difference between panels
in Figure \ref{figure2}. 
One advantage of these methods over those that we have presented in
this paper, is that they do not incur the continuous approximation
error that we see in Figure \ref{figure1}(b) and discussed in
Section \ref{sec:diff}. Diffusion approximation methods would not be
appropriate if one wanted to analyse the dynamics of a slowly changing
variable which has a low average copy number. For
instance, in some gene regulatory network models, there are often
two species, whose number sum to 1 in total, which represent whether a
particular gene is ``switched on'' or ``switched off''. Such a
variable is clearly not a candidate for diffusion
approximation. However, diffusion approximations have been used
successfully for other variables in such
systems~\cite{Kepler:2001:STR}. The dynamics of that gene
cannot themselves be approximated by a diffusion, but may be inferred
from the state of other variables, which may be good candidates for
such an approximation.

One big advantage of the diffusion approximation methods, is the ease
with which the computational effort can be efficiently
parallelised. This is a different approach to parallelisation than in
the case of methods which calculate SSA trajectories. Several trajectories can be computed on different processors,
but in order for the computed invariant distribution to be converged,
either the initial positions of all of the trajectories must be
sampled from the invariant distribution (which is unknown), or the
trajectories must be run for enough time that each one is converged in
distribution. On the other hand, the diffusion approximation
approaches are almost embarrassingly parallelisable, with a subset of
the states for which we wish to estimate the effective drift and
diffusion being given to each processor. The solution of the
Fokker-Planck equation is similarly parallelisable. This means that
given enough computational resources, these algorithms can give us
answers in a very short amount of time.
Moreover, these approaches also allow us to consider adaptive mesh
regimes~\cite{Cotter:2013:AFE}, meaning that one can minimise the
number of sites at which we are required to estimate the effective
drift and diffusion values, while also controlling the global error
incurred. This flexibility is not available in an SSA-type approach.

\end{document}